# DISCUSSION OF "EQUI-ENERGY SAMPLER" BY KOU, ZHOU AND WONG

By Yves F. Atchadé and Jun S. Liu

*University of Ottawa and Harvard University*

We congratulate Samuel Kou, Qing Zhou and Wing Wong (referred to subsequently as KZW) for this beautifully written paper, which opens a new direction in Monte Carlo computation. This discussion has two parts. First, we describe a very closely related method, multicanonical sampling (MCS), and report a simulation example that compares the equi-energy (EE) sampler with MCS. Overall, we found the two algorithms to be of comparable efficiency for the simulation problem considered. In the second part, we develop some additional convergence results for the EE sampler.

**1. A multicanonical sampling algorithm.** Here, we take on KZW's discussion about the comparison of the EE sampler and MCS. We compare the EE sampler with a general state-space extension of MCS proposed by Atchadé and Liu [1]. We compare the two algorithms on the multimodal distribution discussed by KZW in Section 3.4.

Let $(\mathcal{X}, \mathcal{B}, \lambda)$ be the state space equipped with its $\sigma$-algebra and appropriate measure, and let $\pi(x) \propto e^{-h(x)}$ be the density of interest. Following the notation of KZW, we let $H_0 < H_1 < \cdots < H_{K_e} < H_{K_e+1} = \infty$ be a sequence of energy levels and let $D_j = \{x \in \mathcal{X} : h(x) \in [H_j, H_{j+1})\}$, $0 \le j \le K_e$, be the energy rings. For $x \in \mathcal{X}$, define $I(x) = j$ if $x \in D_j$. Let $1 = T_0 < T_1 < \cdots < T_{K_t}$ be a sequence of "temperatures." We use the notation $k^{(i)}(x) = e^{-h(x)/T_i}$, so that $\pi^{(i)}(x) = k^{(i)}(x)/\int k^{(i)}(x)\lambda(dx)$. Clearly, $\pi^{(0)} = \pi$. We find it more convenient to use the notation $\pi^{(i)}$ instead of $\pi_i$ as in KZW. Also note that we did not flatten $\pi^{(i)}$ as KZW did.

The goal of our MCS method is to generate a Markov chain on the space $\mathcal{X} \times \{0, 1, \ldots, K_t\}$ with invariant distribution

$$\pi(x, i) \propto \sum_{j=0}^{K_e} \frac{k^{(i)}(x)}{Z_{i,j}} \mathbf{1}_{D_j}(x) \lambda(dx),$$

where $Z_{i,j} = \int k^{(i)}(x) \mathbf{1}_{D_j}(x) \lambda(dx)$. With a well-chosen temperature sequence $(T_i)$ and energy levels $(H_j)$, such a Markov chain would move very easily from any temperature level $\mathcal{X} \times \{i\}$ to another. And inside each temperature level $\mathcal{X} \times \{i\}$, the algorithm would move very easily from any energy ring $D_j$ to another. Unfortunately, the normalizing constants $Z_{i,j}$ are not known. They are estimated as part of the algorithm using the Wang–Landau recursion that we describe below. To give the details, we need a proposal kernel $Q_i(x, dy) = q_i(x, dy)\lambda(dy)$ on $\mathcal{X}$, a proposal kernel $\Delta(i,j)$ on $\{0, \ldots, K_t\}$ and $(\gamma_n)$, a sequence of positive numbers. We discuss the choice of these parameters later.

ALGORITHM 1.1 (*Multicanonical sampling*).

*Initialization.* Start the algorithm with some arbitrary $(X_0, t_0) \in \mathcal{X} \times \{0, 1, \ldots, K_t\}$. For $i = 0, \ldots, K_t$, $j = 0, \ldots, K_e$ we set all the weights to $\phi_0^{(i)}(j) = 1$.

*Recursion.* Given $(X_n, t_n) = (x, i)$ and $(\phi_n^{(i)}(j))$, flip a $\theta$-coin.

*If Tail.* Sample $Y \sim Q_i(x, \cdot)$. Set $X_{n+1} = Y$ with probability $\alpha^{(i)}(x, Y)$; otherwise set $X_{n+1} = x$, where

$$(1.1) \qquad \alpha^{(i)}(x,y) = \min\left[1, \frac{k^{(i)}(y)}{k^{(i)}(x)} \frac{\phi_n^{(i)}(I(x))}{\phi_n^{(i)}(I(y))} \frac{q_i(y,x)}{q_i(x,y)}\right].$$

Set $t_{n+1} = i$.

*If Head.* Sample $j \sim \Delta(i, \cdot)$. Set $t_{n+1} = j$ with probability $\beta_x(i, j)$; otherwise set $t_{n+1} = i$, where

$$(1.2) \qquad \beta_x(i,j) = \min\left[1, \frac{k^{(j)}(x)}{k^{(i)}(x)} \frac{\phi_n^{(i)}(I(x))}{\phi_n^{(j)}(I(x))} \frac{\Delta(j,i)}{\Delta(i,j)}\right].$$

Set $X_{n+1} = x$.

*Update the weights.* Write $(t_{n+1}, I(X_{n+1})) = (i_0, j_0)$. Set

$$(1.3) \qquad \phi_{n+1}^{(i_0)}(j_0) = \phi_n^{(i_0)}(j_0)(1 + \gamma_n),$$

leaving the other weights unchanged.

If we choose $K_t = 0$ in the algorithm above we obtain the MCS of [9] (the first MCS algorithm is due to [4]) and $K_e = 0$ gives the simulated tempering algorithm of [6]. The recursion (1.3) is where the weights $Z_{i,j}$ are being estimated. Note that the resulting algorithm is no longer Markovian. Under some general assumptions, it is shown in [1] that $\theta_n^{(i)}(j) := \frac{\phi_n^{(i)}(j)}{\sum_{i=0}^{K_t} \sum_{l=0}^{K_e} \phi_n^{(i)}(l)} \to \int_{D_j} \pi^{(i)}(x)\lambda(dx)$ as $n \to \infty$.

The MCS can be seen as a random-scan-Gibbs sampler on the two variables $(x,i) \in \mathcal{X} \times \{0,\ldots,K_t\}$, so the choice $\theta = 1/2$ for coin flipping works well. The proposal kernels $Q_i$ can be chosen as in a standard Metropolis–Hastings algorithm. But one should allow $Q_i$ to make larger proposal moves for larger $i$ (i.e., hotter distributions). The proposal kernel $\Delta$ can be chosen as a random walk on $\{0,\ldots,K_t\}$ (with reflection on 0 and $K_t$). In our simulations, we use $\Delta(0,1) = \Delta(K_t, K_t - 1) = 1$, $\Delta(i, i-1) = \Delta(i, i+1) = 1/2$ for $i \notin \{0, K_t\}$.

It can be shown that the sequence $(\theta_n)$ defined above follows a stochastic approximation with step size $(\gamma_n)$. So choosing $(\gamma_n)$ is the same problem as choosing a step-size sequence in a stochastic approximation algorithm. We follow the new method proposed by Wang and Landau [9] where $(\gamma_n)$ is selected adaptively. Wang and Landau's idea is to monitor the convergence of the algorithm and adapt the step size accordingly. We start with some initial value $\gamma_0$ and $(\gamma_n)$ is defined by $\gamma_n = (1 + \gamma_0)^{1/(k+1)} - 1$ for $\tau_k < n \leq \tau_{k+1}$, where $0 = \tau_0 < \tau_1 < \cdots$ is a sequence of stopping times. Assuming $\tau_i$ finite, $\tau_{i+1}$ is the next time $k > \tau_i$ where the occupation measures (obtained from time $\tau_i + 1$ on) of all the energy rings in all the temperature levels are approximately equal. Various rules can be used to check that the occupation measures are approximately equal. Following [9], we check that the smallest occupation measure obtained is greater than $c$ times the mean occupation, where $c$ is some constant (e.g., $c = 0.2$) that depends on the complexity of the sampling problem.

It is an interesting question to know whether this method of choosing the step-size sequence can be extended to more general stochastic approximation algorithms. A theoretical justification of the efficiency of the method is also an open question.

**2. Comparison of EE sampler and MCS.** To use MCS to estimate integrals of interest such as $\mathbb{E}_{\pi_0}(g(X))$, one can proceed as KZW did by writing $\mathbb{E}_{\pi_0}(g(X)) = \sum_{j=0}^{K_e} p_j \mathbb{E}_{\pi_0}(g(X)|X \in D_j)$. Samples from the high-temperature chains can be used to estimate the integrals $\mathbb{E}_{\pi_0}(g(X)|X \in D_j)$ by importance reweighting in the same way as KZW did. In the case of MCS, the probabilities $p_j = \Pr_{\pi_0}(X \in D_j)$ are estimated by $\hat{p}_j = \frac{\phi_n^{(0)}(j)}{\sum_{l=0}^{K_e} \phi_n^{(0)}(l)}$.

We compared the performances of the EE sampler and the MCS described above for the multimodal example in Section 3.4 of KZW. To make the two samplers comparable, each chain in the EE sampler was run for $N$ iterations. We did the simulations for $N = 10^4$, $N = 5 \times 10^4$ and $N = 10 \times 10^4$. For the

TABLE 1
*Improvement of MCS over EE as given by*
$(\hat{\sigma}_{\mathrm{EE}}(g) - \hat{\sigma}_{\mathrm{MC}}(g))/\hat{\sigma}_{\mathrm{MC}}(g) \times 100$

|  | $\mathbb{E}(X_1)$ | $\mathbb{E}(X_2)$ | $\mathbb{E}(X_2^2)$ | $\mathbf{Pr}_\pi(X \in B)$ |
|---|---|---|---|---|
| $N = 10^4$ | 13.77 | 12.53 | 8.98 | $-63.49$ |
| $N = 5 \times 10^4$ | 6.99 | $-7.31$ | $-10.15$ | $-51.22$ |
| $N = 10^5$ | 1.92 | 5.79 | 4.99 | $-55.31$ |

*The comparisons are based on 100 replications of the samplers for each $N$.

MC sampler, we used $K_t = K_e = K$ and the algorithm was run for $(K + 1) \times N$ total iterations. We repeated each sampler for $n = 100$ iterations in order to estimate the finite sample standard deviations of the estimates they provided. Table 1 gives the improvements (in percentage) of MCS over EE sampling. $\Pr_\pi(X \in B)$ is the probability under $\pi$ of the union of all the discs with centers $\mu_i$ (the means of the mixture) and radius $\sigma/2$. As we can see, when estimating global functions such as moments of the distribution, the two samplers have about the same accuracy with a slight advantage for MCS. But the EE sampler outperformed MCS when estimating $\Pr_\pi(X \in B)$. The MCS is an importance sampling algorithm with a stationary distribution that is more widespread than $\pi$. This may account for the better performance obtained by the EE sampler on $\Pr_\pi(X \in B)$. More thorough empirical and theoretical analyses are apparently required to reach any firmer conclusions.

**3. Ergodicity of the equi-energy sampler.** In this section we take a more technical look at the EE algorithm and derive some ergodicity results. First, we would like to mention that in the proof of Theorem 2, it is not clear to us how KZW derive the convergence in (5). Equation (5) implicitly uses some form of convergence of the distribution of $X_n^{(i+1)}$ to $\pi^{(i+1)}$ as $n \to \infty$ and it is not clear to us how that follows from the assumption that $\Pr(X_{n+1}^{(i+1)} \in A | X_n^{(i+1)} = x) \to S^{(i+1)}(x, A)$ as $n \to \infty$ for all $x$, all $A$.

In the analysis below we fix that problem, but under a more stringent assumption. To state our result, let $(\mathcal{X}, \mathcal{B})$ be the state space of each of the equi-energy chains. If $P_1$ and $P_2$ are two transition kernels on $\mathcal{X}$, the product $P_1 P_2$ is also a transition kernel defined as $P_1 P_2(x, A) = \int P_1(x, dy) P_2(y, A)$. Recursively, we define $P_1^n$ as $P_1^1 = P_1$ and $P_1^n = P_1^{n-1} P_1$. If $f$ is a measurable real-valued function on $\mathcal{X}$ and $\mu$ is a measure on $\mathcal{X}$, we denote $Pf(x) := \int P(x, dy) f(y)$ and $\mu(f) := \int \mu(dx) f(x)$. Also, for $c \in (0, \infty)$ we write $|f| \leq c$ to mean $|f(x)| \leq c$ for all $x \in \mathcal{X}$. We define the following distance between

$P_1$ and $P_2$:

$$\|\|P_1 - P_2\|\| := \sup_{x \in \mathcal{X}} \sup_{|f| \leq 1} |P_1 f(x) - P_2 f(x)|, \tag{3.1}$$

where the supremum is taken over all $x \in \mathcal{X}$ and over all measurable functions $f : \mathcal{X} \to \mathbb{R}$ with $|f| \leq 1$. We say that the transition kernel $P$ is uniformly geometrically ergodic if there exists $\rho \in (0, 1)$ such that

$$\|\|P^n - \pi\|\| = O(\rho^n). \tag{3.2}$$

It is well known that (3.2) holds if and only if there exist $\varepsilon > 0$, a nontrivial probability measure $\nu$ and an integer $m \geq 1$ such that the so-called $M(m, \varepsilon, \nu)$ minorization condition holds, that is, $P^m(x, A) \geq \varepsilon \nu(A)$ for all $x \in \mathcal{X}$ and $A \in \mathcal{B}$ (see, e.g., [8], Proposition 2). We recall that $T_{\mathrm{MH}}^{(i)}$ denotes the Metropolis–Hastings kernel in the EE sampler. The following result is true for the EE sampler.

THEOREM 3.1. *Assume that $\forall i \in \{0, \ldots, K\}$, $T_{\mathrm{MH}}^{(i)}$ satisfies a $M(1, \varepsilon_i, \pi^{(i)})$ minorization condition and that condition* (iii) *of Theorem* 2 *of the paper holds. Then for any bounded measurable function $f$, as $n \to \infty$,*

$$\mathbb{E}[f(X_n^{(i)})] \longrightarrow \pi^{(i)}(f) \quad and \quad \frac{1}{n} \sum_{k=1}^{n} f(X_k^{(i)}) \xrightarrow{a.s.} \pi^{(i)}(f). \tag{3.3}$$

For example, if $\mathcal{X}$ is a compact space and $e^{-h(x)}$ remains bounded away from 0 and $\infty$, then (3.3) holds. Note that the $i$th chain in the EE sampler is actually a nonhomogeneous Markov chain with transition kernels $K_0^{(i)}, K_1^{(i)}, \ldots$, where $K_n^{(i)}(x, A) := \Pr[X_{n+1}^{(i)} \in A | X_n^{(i)} = x]$. As pointed out by KZW, for any $x \in \mathcal{X}$ and $A \in \mathcal{B}$, $K_n^{(i)}(x, A) \to S^{(i)}(x, A)$ as $n \to \infty$, where $S^{(i)}$ is the limit transition kernel in the EE sampler. This setup brings to mind the following convergence result for nonhomogeneous Markov chains (see [5], Theorem V.4.5):

THEOREM 3.2. *Let $P, P_0, P_1, \ldots$ be a sequence of transition kernels on $(\mathcal{X}, \mathcal{B})$ such that $\|\|P_n - P\|\| \to 0$ and $P$ is uniformly geometrically ergodic with invariant distribution $\pi$. Then the Markov chain with transition kernels $(P_i)$ is strongly ergodic; that is, for any initial distribution $\mu$,*

$$\|\|\mu P_0 P_1 \cdots P_n - \pi\|\| \to 0 \quad as\ n \to \infty. \tag{3.4}$$

The difficulty in applying this theorem to the EE sampler is that we do not have $\|\|K_n^{(i)} - S^{(i)}\|\| \to 0$ but only a setwise convergence $|K_n^{(i)}(x, A) - S^{(i)}(x, A)| \to 0$ for each $x \in \mathcal{X}$, $A \in \mathcal{B}$. The solution we propose is to extend Theorem 3.2 as follows.

THEOREM 3.3. *Let $P, P_0, P_1, \ldots$ be a sequence of transition kernels on $(\mathcal{X}, \mathcal{B})$ such that:*

(i) *For any $x \in \mathcal{X}$ and $A \in \mathcal{B}$, $P_n(x, A) \to P(x, A)$ as $n \to \infty$.*

(ii) *$P$ has invariant distribution $\pi$ and $P_n$ has invariant distribution $\pi_n$. There exists $\rho \in (0,1)$ such that $\|\!|P^k - \pi|\!\| = O(\rho^k)$ and $\|\!|P_n^k - \pi_n|\!\| = O(\rho^k)$.*

(iii) *$\|\!|P_n - P_{n-1}|\!\| \leq O(n^{-\lambda})$ for some $\lambda > 0$.*

*Then, if $(X_n)$ is an $\mathcal{X}$-valued Markov chain with initial distribution $\mu$ and transition kernels $(P_n)$, for any bounded measurable function $f$ we have*

$$(3.5) \quad \mathbb{E}[f(X_n)] \longrightarrow \pi(f) \quad \text{and} \quad \frac{1}{n}\sum_{k=1}^n f(X_k) \xrightarrow{a.s.} \pi(f) \qquad \text{as } n \to \infty.$$

We believe that this result can be extended to the more general class of $V$-geometrically ergodic transition kernels and then one could weaken the uniform minorization assumption on $T_{\text{MH}}^{(i)}$ in Theorem 3.1. But the proof will certainly be more technical. We now proceed to the proofs of the theorems. We first prove Theorem 3.3 and use it to prove Theorem 3.1.

PROOF OF THEOREM 3.3. It can be easily shown from (ii) that $\|\!|\pi_n - \pi_{n-1}|\!\| \leq \frac{1}{1-\rho}\|\!|P_n - P_{n-1}|\!\|$. Therefore, Theorems 3.1 and 3.2 of [2] apply and assert that for any bounded measurable function $f$, $\mathbb{E}[f(X_n) - \pi_n(f)] \to 0$ and $\frac{1}{n}\sum_{k=1}^n [f(X_k) - \pi_k(f)] \xrightarrow{a.s.} 0$ as $n \to \infty$. To finish, we need to prove that $\pi_n(f) \to \pi(f)$ as $n \to \infty$. To this end, we need the following technical lemma proved in [7], Chapter 11, Proposition 18.

LEMMA 3.1. *Let $(f_n)$ be a sequence of measurable functions and let $\mu, \mu_1, \ldots$ be a sequence of probability measures such that $|f_n| \leq 1$ and $f_n \to f$ pointwise $[f_n(x) \to f(x)$ for all $x \in \mathcal{X}]$ and $\mu_n \to \mu$ setwise $[\mu_n(A) \to \mu(A)$ for all $A \in \mathcal{B}]$. Then $\int f_n(x)\mu_n(dx) \to \int f(x)\mu(dx)$.*

Here is how to prove that $\pi_n(f) \to \pi(f)$ as $n \to \infty$. By (i), we have $P_n f(x) \to Pf(x)$ for all $x \in \mathcal{X}$. Then, by (i) and Lemma 3.1, $P_n^2 f(x) = P_n(P_n f)(x) \to P^2 f(x)$ as $n \to \infty$. By recursion, for any $x \in \mathcal{X}$ and $k \geq 1$, $P_n^k f(x) \to P^k f(x)$ as $n \to \infty$. Now, write

$$
\begin{aligned}
|\pi_n(f) - \pi(f)| &\leq |\pi_n(f) - P_n^k f(x)| + |P_n^k f(x) - P^k f(x)| \\
&\quad + |P^k f(x) - \pi(f)| \\
&\leq 2\rho^k \sup_{x \in \mathcal{X}} |f(x)| + |P_n^k f(x) - P^k f(x)| \qquad \text{[by (ii)]}.
\end{aligned}
$$
(3.6)

Since $|P_n^k f(x) - P^k f(x)| \to 0$, we see that $|\pi_n(f) - \pi(f)| \to 0$. □

PROOF OF THEOREM 3.1. Let $(\Omega, \mathcal{F}, \mathbb{P})$ be the probability triplet on which the equi-energy process is defined and let $\mathbb{E}$ be its expectation operator. The result is clearly true by assumption for $i = K$. Assuming that it is true for the $(i+1)$st chain, we will prove it for the $i$th chain.

The random process $(X_n^{(i)})$ is a nonhomogeneous Markov chain with transition kernel $K_n^{(i)}(x, A) := \Pr[X_{n+1}^{(i)} \in A | X_n^{(i)} = x]$. For any bounded measurable function $f$, $K_n^{(i)}$ operates on $f$ as follows:

$$K_n^{(i)} f(x) = (1 - p_{\text{ee}}) T_{\text{MH}}^{(i)} f(x) + p_{\text{ee}} \mathbb{E}[R_n^{(i)} f(x)],$$

where $R_n^{(i)} f(x)$ is a ratio of empirical sums of the $(i+1)$st chain of the form

$$\begin{aligned}
R_n^{(i)} f(x) &= \frac{\sum_{k=-N}^{n} \mathbf{1}_{D_{I(x)}}(X_k^{(i+1)}) \alpha^{(i)}(x, X_k^{(i+1)}) f(X_k^{(i+1)})}{\sum_{k=-N}^{n} \mathbf{1}_{D_{I(x)}}(X_k^{(i+1)})} \\
&\quad + f(x) \frac{\sum_{k=-N}^{n} \mathbf{1}_{D_{I(x)}}(X_k^{(i+1)})(1 - \alpha^{(i)}(x, X_k^{(i+1)}))}{\sum_{k=-N}^{n} \mathbf{1}_{D_{I(x)}}(X_k^{(i+1)})}
\end{aligned}$$
(3.7)

[take $R_n^{(i)} f(x) = 0$ and $p_{\text{ee}} = 0$ when $\sum_{k=-N}^{n} \mathbf{1}_{D_{I(x)}}(X_k^{(i+1)}) = 0$], where $\alpha^{(i)}(x, y)$ is the acceptance probability $\min[1, \frac{\pi^{(i)}(y) \pi^{(i+1)}(x)}{\pi^{(i)}(x) \pi^{(i+1)}(y)}]$. $N$ is how long the $(i+1)$st chain has been running before the $i$th chain started. Because (3.3) is assumed true for the $(i+1)$st chain and condition (iii) of Theorem 2 of the paper holds, we can assume in the sequel that $\sum_{k=-N}^{n} \mathbf{1}_{D_{I(x)}}(X_k^{(i+1)}) \geq 1$ for all $n \geq 1$. We prove the theorem through a series of lemmas.

LEMMA 3.2. *For the EE sampler, assumption* (i) *of Theorem 3.3 holds true.*

PROOF. Because (3.3) is assumed true for the $(i+1)$st chain, the strong law of large numbers and Lebesgue's dominated convergence theorem apply to $R_n^{(i)} f(x)$ and assert that for all $x \in \mathcal{X}$ and $A \in \mathcal{B}$, $K_n^{(i)}(x, A) \to S^{(i)}(x, A)$ as $n \to \infty$, where $S^{(i)}(x, A) = (1 - p_{\text{ee}}) T_{\text{MH}}^{(i)}(x, A) + p_{\text{ee}} \sum_{j=0}^{K} T_{\text{EE}}^{(i,j)}(x, A) \mathbf{1}_{D_j}(x)$, where $T_{\text{EE}}^{(i,j)}$ is the transition kernel of the Metropolis–Hastings with proposal distribution

$$\pi^{(i+1)}(y) \mathbf{1}_{D_j}(y) / p_j^{(i+1)}$$

and invariant distribution

$$\pi^{(i)}(x) \mathbf{1}_{D_j}(x) / p_j^{(i)}. \qquad \square$$

LEMMA 3.3. *For the EE sampler, assumption* (ii) *of Theorem 3.3 holds.*

PROOF. Clearly, the minorization condition on $T_{\mathrm{MH}}^{(i)}$ transfers to $K_n^{(i)}$. It then follows that each $K_n^{(i)}$ admits an invariant distribution $\pi_n^{(i)}$ and is uniformly geometrically ergodic toward $\pi_n^{(i)}$ with a rate $\rho_i = 1 - (1 - p_{\mathrm{ee}})\varepsilon_i$. The limit transition kernel $S^{(i)}$ in the EE sampler as detailed above has invariant distribution $\pi^{(i)}$ and also inherits the minorization condition on $T_{\mathrm{MH}}^{(i)}$. □

LEMMA 3.4. *For the EE sampler, assumption* (iii) *of Theorem* 3.3 *holds true with* $\lambda = 1$.

PROOF. Any sequence $(x_n)$ of the form $x_n = \frac{\sum_{k=1}^n \alpha_n u_n}{\sum_{k=1}^n \alpha_n}$ can always be written recursively as $x_n = x_{n-1} + \frac{\alpha_n}{\sum_{k=1}^n \alpha_k}(u_n - x_{n-1})$. Using this, we easily have the bound

$$|K_n^{(i)} f(x) - K_{n-1}^{(i)} f(x)| \leq 2\mathbb{E}\left[\frac{1}{\sum_{k=-N}^n \mathbf{1}_{D_{I(x)}}(X_k^{(i+1)})}\right]$$

for all $x \in \mathcal{X}$ and $|f| \leq 1$. Therefore, the lemma will be proved if we can show that

$$(3.8) \qquad \sup_{0 \leq j \leq K} \mathbb{E}\left[\frac{n}{\sum_{k=-N}^n \mathbf{1}_{D_j}(X_k^{(i+1)})}\right] = O(1).$$

To do so, we fix $j \in \{0, \ldots, K\}$ and take $\varepsilon \in (0, \delta)$, where $\delta = (1 - p_{\mathrm{ee}})\varepsilon_{i+1} \times \pi^{(i+1)}(D_j) > 0$. We have

$$\mathbb{E}\left[\frac{n}{\sum_{k=-N}^n \mathbf{1}_{D_j}(X_k^{(i+1)})}\right]$$

$$(3.9) \quad = \mathbb{E}\left[\frac{n}{\sum_{k=-N}^n \mathbf{1}_{D_j}(X_k^{(i+1)})} \mathbf{1}_{\{\sum_{k=-N}^n \mathbf{1}_{D_j}(X_k^{(i+1)}) > n(\delta - \varepsilon)\}}\right]$$

$$+ \mathbb{E}\left[\frac{n}{\sum_{k=-N}^n \mathbf{1}_{D_j}(X_k^{(i+1)})} \mathbf{1}_{\{\sum_{k=-N}^n \mathbf{1}_{D_j}(X_k^{(i+1)}) \leq n(\delta - \varepsilon)\}}\right].$$

The first term on the right-hand side of (3.9) is bounded by $1/(\delta - \varepsilon)$.

The second term is bounded by

$$n \Pr\left[\sum_{k=-N}^0 \mathbf{1}_{D_j}(X_k^{(i+1)}) + \sum_{k=1}^n (\mathbf{1}_{D_j}(X_k^{(i+1)}) - \delta) \leq -n\varepsilon\right]$$

$$(3.10) \qquad \leq n \Pr[M_n^{(i+1)} \geq n\varepsilon],$$

DISCUSSION 9

where $M_n^{(i+1)} = \sum_{k=1}^n K_{k-1}^{(i+1)}(X_{k-1}^{(i+1)}, D_j) - \mathbf{1}_{D_j}(X_k^{(i+1)})$. For the inequality in (3.10), we use the minorization condition $K_{k-1}^{(i+1)}(x, D_j) \geq \delta$. Now, the sequence $(M_n^{(i+1)})$ is a martingale with increments bounded by 1. By Azuma's inequality ([3], Lemma 1), we have $n \Pr[M_n^{(i+1)} \geq n\varepsilon] \leq n \exp(-n\varepsilon^2/2) \to 0$ as $n \to \infty$. $\square$

Theorem 3.1 now follows from Theorem 3.3. $\square$
## REFERENCES

[1] ATCHADÉ, Y. F. and LIU, J. S. (2004). The Wang–Landau algorithm for Monte Carlo computation in general state spaces. Technical report.
[2] ATCHADÉ, Y. F. and ROSENTHAL, J. S. (2005). On adaptive Markov chain Monte Carlo algorithms. *Bernoulli* **11** 815–828. MR2172842
[3] AZUMA, K. (1967). Weighted sums of certain dependent random variables. *Tôhoku Math. J.* (*2*) **19** 357–367. MR0221571
[4] BERG, B. A. and NEUHAUS, T. (1992). Multicanonical ensemble: A new approach to simulate first-order phase transitions. *Phys. Rev. Lett.* **68** 9–12.
[5] ISAACSON, D. L. and MADSEN, R. W. (1976). *Markov Chains*: *Theory and Applications*. Wiley, New York. MR0407991
[6] MARINARI, E. and PARISI, G. (1992). Simulated tempering: A new Monte Carlo scheme. *Europhys. Lett.* **19** 451–458.
[7] ROYDEN, H. L. (1963). *Real Analysis*. Collier-Macmillan, London. MR0151555
[8] TIERNEY, L. (1994). Markov chains for exploring posterior distributions (with discussion). *Ann. Statist.* **22** 1701–1762. MR1329166
[9] WANG, F. and LANDAU, D. P. (2001). Efficient, multiple-range random walk algorithm to calculate the density of states. *Phys. Rev. Lett.* **86** 2050–2053.

DEPARTMENT OF MATHEMATICS
AND STATISTICS
UNIVERSITY OF OTTAWA
585 KING EDWARD STREET
OTTAWA, ONTARIO
CANADA K1N 6N5
E-MAIL: yatchade@uottawa.ca

DEPARTMENT OF STATISTICS
SCIENCE CENTER
HARVARD UNIVERSITY
CAMBRIDGE, MASSACHUSETTS 02138
USA
E-MAIL: jliu@stat.harvard.edu